\documentclass[12pt]{article}
\usepackage[margin=1in]{geometry}
\usepackage{amsmath,amsthm,amssymb}

\usepackage{graphicx}
\usepackage{float}
\usepackage{blindtext}

\theoremstyle{plain}
\newtheorem{thm}{Theorem}

\newtheorem{lemma}[thm]{Lemma}

\theoremstyle{definition}

\theoremstyle{remark}
\newtheorem{rem}[thm]{Remark}
\newtheorem*{ex}{Example}

\numberwithin{equation}{section}
\numberwithin{thm}{section}

\usepackage{wasysym} 

\title{Higher Order Bipartiteness vs Bi-Partitioning in Simplicial Complexes}

\author{Marzieh Eidi$^{1,2,*}$, Sayan Mukherjee$^{1,2,3,\dagger}$ \\ 
$^1$Center for Scalable Data Analytics and Artificial Intelligence, Leipzig University\\
$^2$Max Planck Institute for Mathematics in the Sciences\\
$^3$Duke University \\ $*$ meidi@mis.mpg.de\\
$\dagger$ sayan.mukherjee@mis.mpg.de}

\begin{document}

\maketitle

\begin{abstract}
Bipartite graphs are a fundamental concept in graph theory with diverse applications. A graph is bipartite iff it contains no odd cycles, a characteristic that has many implications in diverse fields ranging from matching problems to the construction of complex networks.  Another key identifying feature is  their Laplacian spectrum as bipartite graphs achieve the maximum possible eigenvalue of graph Laplacian.  However, for modeling higher-order connections in complex systems, hypergraphs and simplicial complexes are required due to the limitations of graphs in representing pairwise interactions. In this article, using simple tools from graph theory, we extend the cycle-based characterization from bipartite graphs to those simplicial complexes  that achieve the maximum  Hodge Laplacian eigenvalue, known as disorientable simplicial complexes. We show that a $N$-dimensional simplicial complex is disorientable if its down dual graph contains no simple odd cycle of distinct edges and no twisted even cycle of distinct edges. 
 Furthermore, we see that  in a $N$-simplicial complex without twisting cycles, 
   the fewer the number of (non-branching) simple odd cycles in its down dual graph, the closer is its maximum eigenvalue to the possible maximum eigenvalue of Hodge Laplacian. Similar to the graph case, the absence of odd cycles plays a crucial role in solving the bi-partitioning problem of simplexes in higher dimensions.

\end{abstract}
\section{introduction}

.

Bipartite graphs are useful tools in diverse domains, from matching problems and coding theory to social networks and biomedical applications such as cancer detection\cite{book,cancer,medicine}. 
Exploring the properties of bipartite graphs, and their applications not only has enriched graph theory but also has opened the door to innovative solutions in challenging applications \cite{book}. 
In the past few years, there has been a trend in understanding cycles in networks in terms of function, dynamics, and synchronizability \cite{cycle1}. Cycles and bipartite graphs are closely related. 
A graph is bipartite if and only if it has no odd cycles. Cycle-based characterization of bipartite graphs  is essential in various domains, such as clustering, coloring and matching problems as well as complex network construction and  analysis where detecting quasi-bipartite clusters, i.e. graphs with a few number of even cycles, is needed \cite{PhysRevE1, PhysRevE2}. Moreover, this characterization as well as some of its applications have been extended to directed graphs  where a strongly connected digraph is bipartite if and only if it
has no cycle of odd length \cite{digraph}.

One the other side, one can test for a graph being bipartite
by observing the spectrum of the graph Laplacian, a matrix computed based on the incidence relations between vertices and edges of the graph. 
For the normalized Laplacian, the graph is bipartite if and only if two is the maximum eigenvalue of its Laplacian. Consider the spectrum of the normalized Laplacian of an unweighted graph in increasing order, $ 0 =\lambda_1 \leq...\leq \lambda_n\leq 2$; the graph is connected if and only if  $\lambda_2$ is bigger than zero and is bipartite iff  $\lambda_n$ is equal to two. While the down side of the spectrum tells the number of connected components (i.e. graph's topology), the up side reveals bipartiteness.  
Both the top and bottom eigenvalues then
can help us to partition the vertices of the graph into two sets in two different ways; if $\lambda_2$ is positive but small we can partition the vertices into two sets with very few edges connecting the sets while elements of each set are highly connected to each other and they form clusters.  And when $\lambda_n$ is close to two we can partition the set of the vertices into two sets in a way that for each of these sets, there are few connections, and almost all of the edges, connect the vertices of one set to the other (Fig.1.).   
\begin{figure}[ht]
\centering
\includegraphics[width=0.70\textwidth]{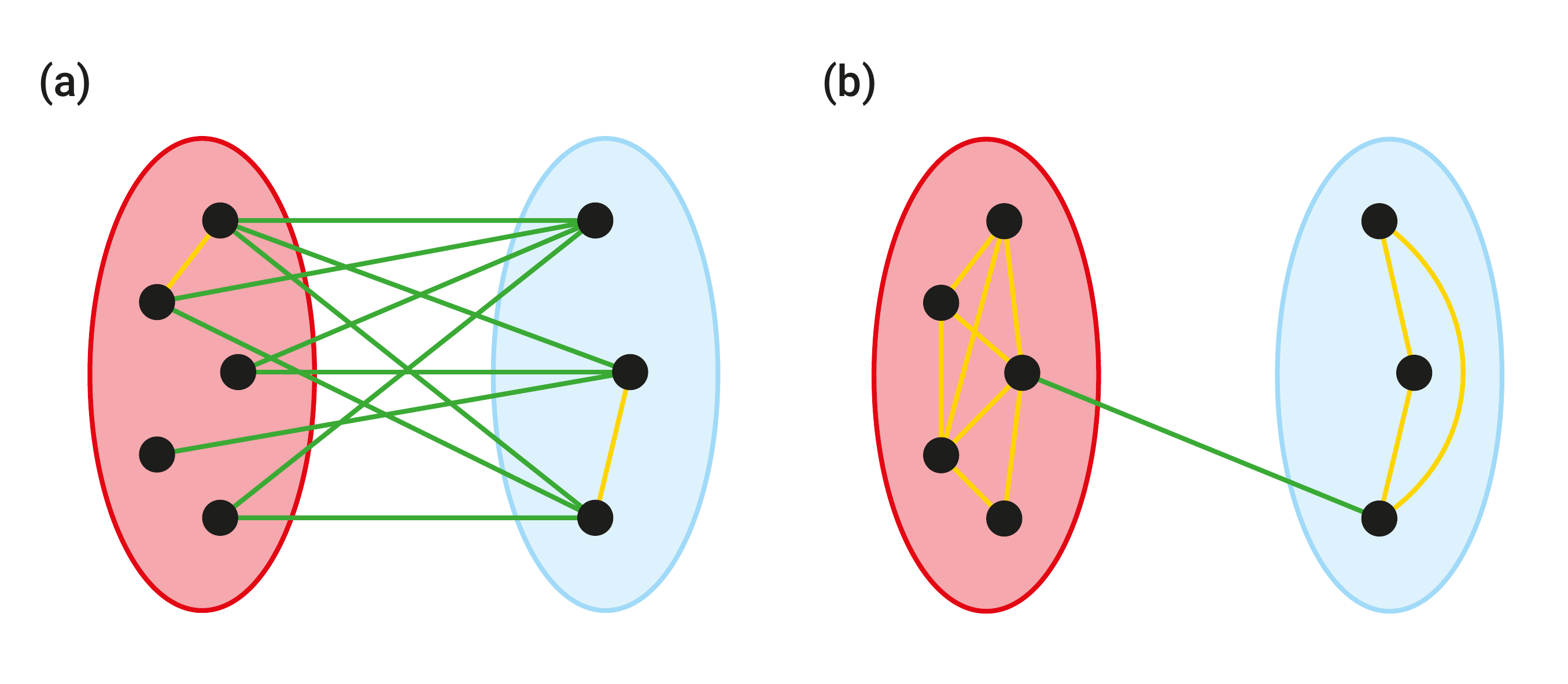}
\caption{Partitioning vertices based on the Laplacian spectrum. In (a) vertices are partitioned based on $\lambda_n$ and in (b) based on the $\lambda_2$.}
\end{figure}

Many empirical networks and complex systems however incorporate higher-order relations between elements and therefore are modeled as hypergraphs and/or simplicial complexes, rather than graphs. Simplicial complexes are generalizations of graphs where there are not only vertices and edges but also possibly triangles, tetrahedrons, and so on.
To develop a systematic tool for the structural analysis of simplicial complexes, different methods, and theories have been extended from graphs to these higher-order structures. 
An important question is that can the graph partitioning scheme that we described, be generalized to simplicial complexes?
The main starting point to answer this question is the discrete "Hodge" Laplacian which is a generalization of the graph Laplacian for simplicial complexes. There have been many advances in our understanding of the spectrum of discrete Hodge Laplacian on simplicial complexes in the past few years \cite{Max, Jost, Horak}. Moreover,  Laplacian-based methods have become popular for detecting the structure and dynamics of complex networks modeled by simplicial complexes in the past few years\cite{Gin}. 
It is well known, due to Eckmann \cite{Eckmann}, that for any fixed $d$ where $N \geq d \geq 0$ ($N$ is the dimension of the complex), the minimum eigenvalue of the $d$-Laplacian can tell us about the topology of the complex in dimension $d$. This theorem has led to many theoretical findings regarding the minimum eigenvalues of the Laplacian on simplicial complexes\cite{Jost, Horak} as well as clustering methods for $d$-simplexes in such a way that the clusters represent the topology of the complex in dimension $d$ \cite{ebli, Schaub_2020}. This is a generalization of what was described at the beginning based on the down side of the Laplacian spectrum for clustering the vertices of a graph. However, there is very little known about the maximum eigenvalue of the Laplacian and specifically the higher dimensional analogues for bipartiteness. A higher-order notion of bipartiteness called "disorientability", was introduced in \cite{Ros} as a structure that the spectrum of the (up) Hodge Laplacian achieves its maximum possible value; the graph is the simplest setting and a one-dimensional simplicial complex is disoriantable if and only if it is bipartite. It is known that simple random walks on the vertices of a connected bipartite graph are periodic. It has been shown that the same is true when the random walks are on higher dimensional simplexes of a general simplicial complex and such random walks are periodic if and only if the complex is disorientable \cite{Ros, Sayan, eidi2023}.  Random walks on graphs and simplicial complexes are the main tools in tackling diverse complex real-world problems from ranking web pages in Google page rank algorithm \cite{Page1} to signal processing and flow network decomposition\cite{Schaub_2020}. 
Considering the numerous theoretical implications and practical applications that bipartite graphs have, it is very desirable to develop our understanding of disorientability in simplicial complexes not just in terms of the spectrum of the Laplacian, but in terms of the combinatorial structure and shape of the complex and the parity of length of higher order cycles and see if we can have similar simplified characterization for recognizing whether a simplicial complex is disorientable and does disorientability give some kind of bi-portioning of the simplicies of higher dimensions in a similar manner like bipartite graphs? Similarly, can we  detect quasi-bipartite clusters in higher dimensions? Since more data sets and some complex problems are nowadays represented by simplicial complexes and not graphs it seems necessary to develop the applied partitioning methods for partitioning the simplexes of dimensions higher than zero and explore further the applications, and devise algorithms.

 In this manuscript, we state conditions under which the simplicial complex is disorientable, in terms of the length of the higher dimensional cycles in the complex. Consequently, any simplicial complex can become disorientable by a finite number of splittings of some of the maximum dimensional simplexes; note that this is very similar to the graph case where by dividing some of the edges into two, and playing with the parity of the length of their cycles,  we can always make the graph bipartite. And the fewer number of splittings needed to make the graph bipartite, the closer the maximum eigenvalue would be to the possible maximum eigenvalue of the graph Laplacian. We elaborate this for graphs and simplicial complexes in more details in the next sections (Lemma 3.8 and Theorem 3.9). \\
To the best of our knowledge, this is the first full characterization of
disorientability from the combinatorial structure of the higher-dimensional cycles in simplicial complexes. We 

First, we recall some preliminary notions. 

\section{Background}
\textit{Simplicial complexes-}
A simplicial complex $K$ on a vertex set
$V=\lbrace v_1, ... , v_n \rbrace$ consists of a collection of simplices, that is, subsets of $V$ with the requirement that all these subsets are closed under inclusion. A subset that contains $d + 1$-vertices, is called a $d$-simplex. A $0$-simplex is simply a
vertex, a $1$-simplex is an edge, and a $2$-simplex is a triangle. The dimension of the complex is the maximum $d$, where we have at least one $d$-simplex in K.
For computational purposes, we need to define an orientation for each d-simplex when $d\geq$ 1. An orientation of a $d$-simplex  is an equivalence class of
orderings of its vertices, where two orderings are equivalent if they differ by an even
permutation. For simplicity, we can choose the reference orientation of the simplices
induced by the ordering of the vertex labels. A $0$-simplex (a node) can have only one orientation. Hence, issues of orientation do not arise in graph-theoretic settings.
 Let $S_d$ be the set of all $d$-simplexes in $K$ with   $0 \leq d \leq N$ and $[S_d]$ be the set of all oriented $d$-simplexes. For any $d$-simplex $\sigma$ ($d\neq0$) we have two opposite orientations, clockwise and counterclockwise (note that both are in $[S_d]$). 
 A face of a $d$-simplex $\sigma$ 
is a subset of $\sigma$
 with cardinality $d$, i.e., with one element
of $\sigma$
 omitted. If $\rho$
is a face of a $d$-simplex $\sigma_d$, then $\sigma_d$ is called a co-face of $\rho_{d-1}$. The degree of a $d$-simplex is the number of its co-faces (i.e., of dimension $d+1$). A $d$-simplex is called branching if its degree is bigger than $2$.
      A $d$-cycle of length $l$ ($0\leq d \leq N-1$) is a chain of $d$ simplexes,  $\sigma^0_d...,\sigma^l_{d}$ such that for each $i$, $\sigma^i$ and $\sigma^{i+1}$ are upper adjacent (i. e. they share a co-face) and $\sigma^0=\sigma^l$.

A  $d$-cycle is called non-twisted (or simple) if there exist an assignment of orientations on its $d+1$-simplexes such that the induced orientation of every two $d+1$-simplex on their common $d$-face cancel each other (namely they are opposite). Otherwise we call it twisted.   We call the simplicial complex non-twisting if it does not have a twisted cycle 
 and otherwise, we call it twisting. A discrete cylinder and torus with respectively simple $0$-cycle(s) and $1$-cycle(s) are examples of non-twisting simplicial complexes and discrete Mobius strip and  Klein bottle are twisting as they have respectively a $0$-twisting cycle and a $1$-twisting cycle (see Fig 4).


      Note that graphs can not have twisted cycles and this just can happen in higher dimensions. Twisting cycles do not naturally arise in geometric and topological data analysis (TDA) when constructing Simplicial complexes from a set of data points. However, we present it here for the sake of completeness of our characterization. \\
\textit{Boundary and co-boundary matrices and the Laplacian-} We can extend face/co-face relations to the oriented $d$-simplexes with the help of the boundary/coboundary operators.   
The $d$-th chain group $C_d(K)$ of $K$ is a vector space with real coefficients with the basis $S_d$. The {\em boundary map } 
$\partial_d: C_d(K) \rightarrow C_{d-1}(K)$ is a linear operator defined by 
\begin{equation} 
\partial_d[i_0,...i_d]=\sum_{j=0}^{d} (-1)^j [i_0,. i_{j-1},i_{j+1}..i_d].
\end{equation}
After choosing a basis for $C_d$, the boundary operator $\partial_d$ can be represented by a matrix $B_d$, which enables us to simply perform computations via matrix calculus. For graphs, the matrix $B_1$ is the node-to-edge incidence matrix. 
Likewise, the higher-order boundary maps can be interpreted
as higher-order incidence matrices between simplices and their (co-)faces for each $B_i$. The transpose of the boundary matrix $B_i$ denoted by $B_i^T$ represents the co-boundary map $\partial_d^T: C_d \rightarrow C_{d+1}$. Subsequently, we can define the $d$-th discrete (combinatorial) Hodge Laplacian  as 
\begin{equation}
L_d=B_d^T B_d + B_{d+1}B_{d+1}^T.
\end{equation}
For $d=0$, $L_0= B_{1}B_{1}^T$ and for $d=N$, $L_N= B_N^T B_N$. $B_d^T B_d$  and $B_{d+1}B_{d+1}^T$ are respectively called the down and up Laplacians. 
\begin{rem}
    For every $d$ the  non-zero spectrum of $L^{up}_d$ equals to the non-zero spectrum of $L^{down}_{d+1}$  \cite{Horak}. 
\end{rem}
The  Laplacian spectrum and in particular its maximum eigenvalue can be obtained using the Rayleigh quotient. For a graph G we have: \\
\begin{center}
    $\lambda_n= Sup_f \frac{\sum_{u \sim v} (f(u)-f(v))^2}{\sum_{v} f^2(v) d(v)} $
\end{center}
where  $\lambda_n$ is the maximum  eigenvalue of the normalized Laplacian and the supremum is taken over  real-valued functions on the vertex set and $d(v)$ is the degree of the vertex $v$. 
This formula can be extended to weighted graphs as well as higher dimensions. For a $N$-dimensional simplicial complex $\lambda_n \leq N+1$ \cite{Jost}.

\textit{Dual graphs-} We can create dual graphs of simplicial complexes based on the upper adjacency and/or lower adjacency connections between the simplexes. For this purpose, we consider a N-simplicial complex $K$ and we fix a number $d$ ( $0 \leq d \leq N)$.
The up-dual graph $G$ of $K$ in dimension $d$ is constructed as follows: 
each $d$-simplex in $K$ becomes a vertex in the dual graph
$G$, and there is an edge between two vertices in $G$ if their corresponding d-simplexes  in $K$ share a co-face. Similarly, the down-dual graph $G$ of $K$ in dimension $d$ is a graph $G$ that its vertices are the $d$-simplexes of $K$ but there is an edge between two vertices in $G$ if their corresponding d-simplexes in $K$ share a face. Particularly, the down-dual graph of a graph is called a line graph. We call a $N$-complex $d$-connected if the up dual graph in dimension $d$ is connected.

\textit{Higher order bipartiteness-}  There are two types of generalizations for bipartiteness for a $N$-simplicial complex; a notion which has been investigated more, as it is connected to random walks \cite{Ros},  is called disorientability. A disorientation of a $N$-complex $K$ is a choice of orientations of its $N$-simplexes, that whenever two arbitrary $N$-simplexes intersect in a $(N-1)$-simplex, they induce the same orientation on it (see Fig.2.). If $K$ has a disorientation it is said to be disorientable. A graph is disorientable iff it is bipartite \cite{Ros}.   
\begin{figure}[ht]
\centering
\includegraphics[width=0.35\textwidth]{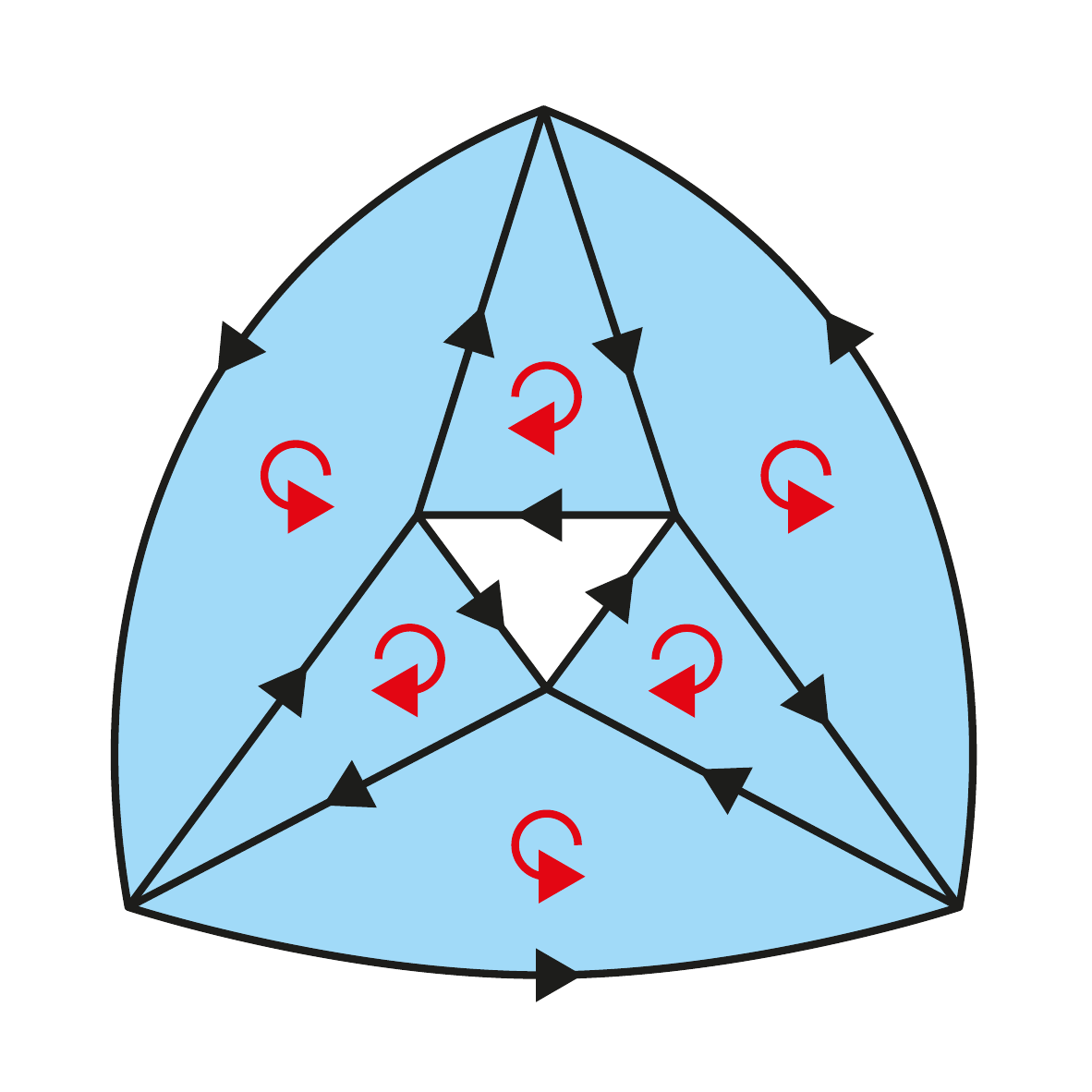}
\caption{Example of a 2-d disorientable simplicial complex. 
}
\end{figure}

Another natural analogue is
“$(d + 1)$-partiteness”: having some partition $A_0, . . . , A_d$ of $V$ so that every $d$-simplex contains
one vertex from each $A_i$. A $(d + 1)$-partite complex is seen to be disorientable, but the opposite does not necessarily hold for $d\geq 2$ \cite{Ros}. Since this second case is not related to the spectrum of the Laplacian and periodicity of random walks and does not have the theoretical implications that disorientability has we will omit it here.

\section{Higher Order Bipartiteness}
To explore disorientable $N$-dimensional simplicial complexes we use their down-dual graphs and we start with graphs as 1-d complexes. As already mentioned in \cite{Ros} bipartite graphs are disorientable and in fact, these two are equivalent in graphs. We elaborate on this from a signs perspective and the line graph of a graph. \\
If a graph is disorientable, there is a choice of orientations on its edges that adjacent edges induce the same orientation on their common vertex. Equivalently,  there is an assignment of $\pm$ to vertices such that no two adjacent vertex have the same sign, namely we have a bipartite graph where each partition is labeled with one of these signs representing the head/tails of the oriented edges.   Now let's look at the line graph of a bipartite graph. We have the following simple observation:

    When going to line graphs of a general graph, its cycles and the parity of their length are preserved. But we might get some more cycles in the line graph that did not exist before due to the existence of branching vertices, i.e. those vertices that have degree bigger or equal than 3. We call these cycles branching cycles as they correspond to the branching vertices in the original graph. For any such vertex $v$ with degree $\alpha \geq 3 $ we would have a sequence of cycles with lengths $\alpha$, $\alpha-1$,..,  3 in the line graph.   So if there exists no branching vertex, there is a one-to-one correspondence between the cycles in the graphs and its dual (line graph). Therefore if the graph is bipartite and non-branching its line graph is also bipartite.

Therefore a non-branching disoientable graph is a graph where we can assign $+$ and $-$ to the edges of its line graph such that no two adjacent edges have the same sign. But what if the bipartite graph has some branching vertices? In this case, the line graph is not bipartite and has (fundamental) odd cycles that for every such cycle,  all of its edges relate to the same corresponding branching vertex. We note that if we consider general (not necessarily fundamental) cycles the existing odd cycles include more than one edge corresponding to that branching vertex. \\
Therefore, the graph is disorientable iff it is bipartite or equivalently iff its line graph has no odd cycle of distinct edges, namely they correspond to different vertices in the original graph. In the language of signs, a graph is disorientable iff in its line graph, we can assign + and - to the edges such that no two adjacent edges that correspond to different vertices in the original graph, get the same sign. 

\begin{figure}[ht]
\centering
\includegraphics[width=0.60\textwidth]{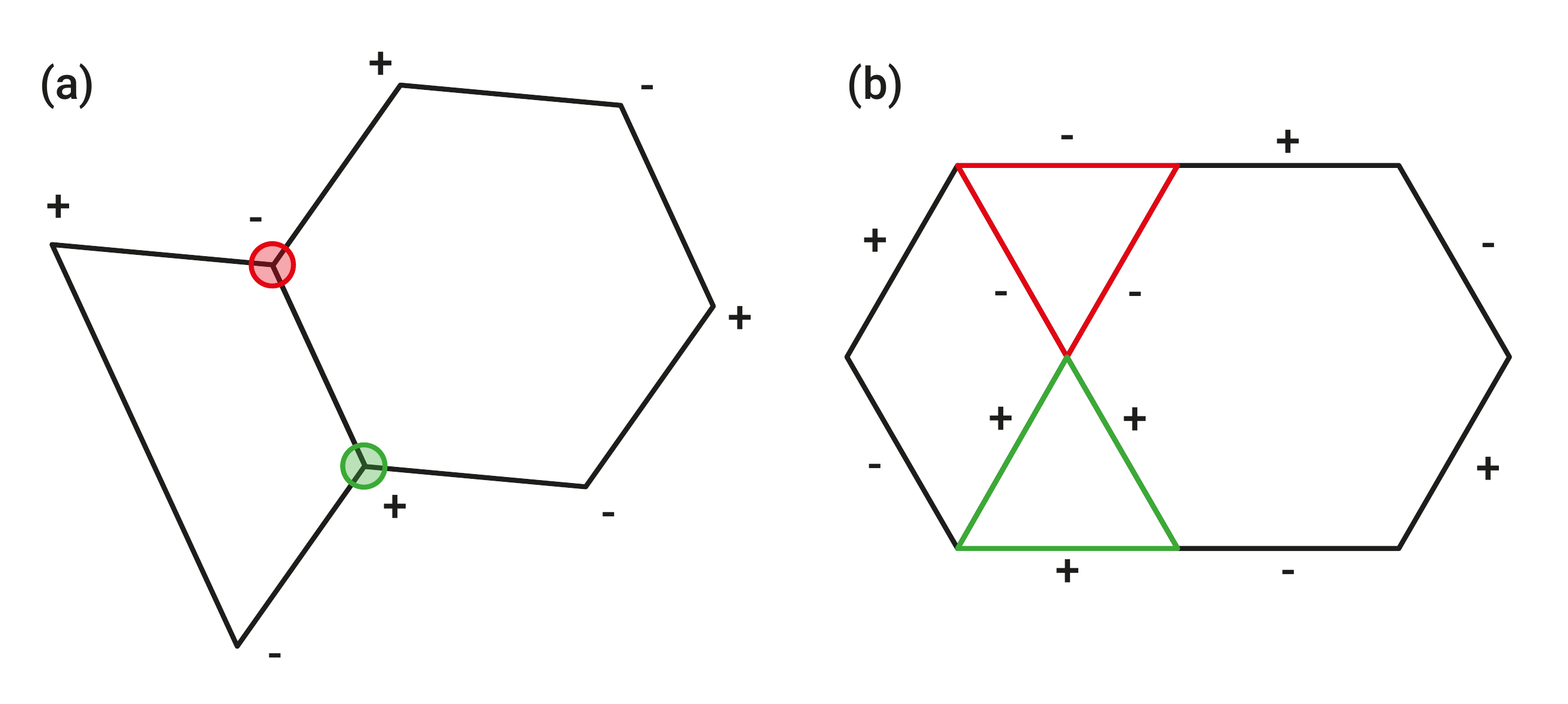}
\caption{A bipartite graph and its line graph. (a) is bipartite and  (b) is the line graph of (a) and has two odd cycles of length 3 corresponding to the branching red and green vertices in (a).}
\end{figure}

This simple but fundamental change of perspective from graphs to their line graphs helps us to explore bipartiteness in general simplicial complexes. We start with the simplest cases: non-branching and non-twisting simplicial complexes and we gradually develop the idea for the general case.

\begin{enumerate}
    \item If the simplicial complex is non-branching and non-twisting: 
    \begin{lemma}
        A non-branching, non-twisting simplicial complex is disorientable if and only if in its down dual graph, we can assign + and - to the edges in such a way that no two adjacent edges have the same sign. We note that this is equivalent to the condition that the down-dual graph has no odd cycle which is the same as bipartiteness.  
    \end{lemma} 
  
\begin{proof}
  If the simplicial complex is disorientable, by definition there is an assignment of orientations on $N$-simplexes such that they induce the same orientation on their common face. Therefore, in its corresponding down dual graph, we would be able to have a single assignment of +/- to each edge representing the induced orientations on the common $(N-1)$-faces (in the simplicial complex) such that no two adjacent edges get the same sign. This means that the down-dual graph has no odd cycle. Therefore it is bipartite.
  For the reverse, if there is such assignment of +/- to the edges of the down-dual graph, by definition the simplicial complex is disorientable as each such sign for each edge can be interpreted as the unique induced sign on the common $(N-1)$-simplex (from its oriented N-simplexes cofaces).
\end{proof}
\begin{rem}
    The above lemma can be considered as a higher-dimensional version of the two-coloring problem on graphs. We note that in this case, we can color the max-dimensional simplexes with two colors in such a way that no two adjacent simplexes that share a face have the same color. Recall that a graph can be colored by two colors if and only if it is bipartite. Also planer graphs (graphs that can be drawn without any of their edges crossing) can be colored using at most four colors, such that no two adjacent nodes have the same color.
\end{rem}
    \item  
The simplicial complex is branching and non-twisting:

\begin{lemma}
    
 A branching non-twisting simplicial complex is disorientable if and only in its down dual graph, either all of the edges of an odd cycle correspond to a branching  $(N-1)$-simple and/or that odd cycle includes more than one edge corresponding to that
branching simplex.
\end{lemma}
\begin{proof}
 Similar to the branching bipartite graphs that are already described, for every branching $(N-1)$-simple $\alpha$, in the down dual graph we would have a (branching) cycle of length= degree $\alpha$, as well as all (branching) cycles with length between three and degree of $\alpha$, such that all of the edges of these cycles correspond to $\alpha$.  Therefore, all of the edges corresponding to all of these cycles in the down dual graph are adjacent and get the same sign and when moving to another adjacent edge (that corresponds to another $(N-1)$-simple), the sign is changed. Therefore, we might have odd cycles, but they can just be for the branching $(N-1)$-simplexes, and no other fundamental odd cycle will exist. Also, when considering general odd cycles (not necessarily fundamental ones) they might include more than one edge corresponding to that
branching simplex.   
    \end{proof}
    \item The simplicial complex is non-branching and twisting:
       \begin{lemma}
        A non-branching simplicial complex is disorientable if and only if in its down dual graph, all the twisting cycles have odd length and all the non-twisting cycles have even length. 
    \end{lemma}
    \begin{proof}
As we saw, in non-branching (disoriantable) simplicial complexes, the down dual graph has no simple cycle of odd length. We should show every twisted cycle in the down dual graph should be odd to assure disorientability. We note that twisted cycles, by definition, always include two $N$-simplexes inducing the same orientation on their common face. We recall that this happens whenever we try to assign orientations on the N-simplexes in such a way that  lower-adjacent $N$-simplexes induce opposite orientation on their common face. We call this three simplexes the twisted part. Lets assume the twisted cycle is as following:\\
  $(N-1)$-cycle of length $l$,  $\sigma^0_{N-1}, \beta^1_{N} ...,\beta^{l}_{N},\sigma^l_{N-1}$ \\ where $\sigma^0=\sigma^l$ and $\sigma^0 ( \sigma^l)$, $\beta^{1}_{N}$ and $\beta^{l}_{N}$ is the twisting part.

 If we start assigning orientations on $N$-simplexes, starting from  $\beta^{1}_{N}$ in a compatible manner where every two adjacent $N$-simplex in the cycle induce the same orientation on their common face, if $l$ is even, due to the twisting part, the induced orientation of $\beta^{l}_{N}$  on $\sigma^l$ would be opposite of orientation of $\sigma^0$ induced by $\beta^{1}_{N}$ (see Fig 4.a). Therefore we need odd number of $N$-simplexes in the cycle to make compatible induced orientations on all common faces along the cycle and in particular on the $\sigma^0= \sigma^l$ (see Fig 4.c).
Roughly speaking, the twisted part acts as an extra (hidden) $N$-dimensional simplex where the twisted connection corresponds to two faces of this simplex, with opposite orientations. Therefore, every twisted cycle as opposed to simple cycles should have odd length to induce the dissorientability condition along the cycle. The rest of the proof is as before and we omit it here.
    \end{proof}

\end{enumerate}
In conclusion, we have the following theorem that enables us to check the dissorientability of a general simplicial complex based on the length of the cycles of its down dual graph.
\begin{thm}
     A simplicial complex is dissorientable if and only if in its down dual graph, the (possible) odd cycles only correspond to the branching $(N-1)$-simplexes and/or twisted cycles and there is no twisted cycle of even length.
\end{thm}
\begin{proof}
    This can be simply obtained by combining Lemma 3.1, Lemma 3.3 and Lemma 3.4.
\end{proof}
\begin{thm}
    Every simplicial complex can become dissorientable by a finite number of splittings of some of the $N$-dimensional simplexes into two $N$-dimensional simplexes.
\end{thm}
 \begin{proof}
       We orient all of $N$-simplexes one by one in a compatible manner, namely in such a way that the adjacent $N$-simplexes induce the same orientation on their common face. If we can achieve this goal globally and for all of the lower-adjacent $N$-simplexes, we are done; i.e., the simplicial complex is disorientable. If not, there are at least two $N$-simplexes, $A$ and $B$, such that they induce opposite orientations on their common face. We note that this means $A$ and $B$ are adjacent vertices in the down dual graph and based on the previous theorem both are included in a cycle of odd (even for the twisting cycles) length. We chose one of them arbitrarily (A) and we split it into two $N$-simplexes $A'$ and $A''$, where $A''$ and $B$ have non-empty intersections. We orient $A'$ based on the orientation of A in such a manner that $A$ and $A'$ have exactly the same orientation on  the shared face. We then orient $A''$ in a compatible manner with $A'$ where they induce the same orientation on their common face. Then trivially $A''$ will have a compatible orientation with $B$, meaning they induce the same orientation on their common face. If we repeat this process for any two non-compatible N-simplex we get our desired result. We note that by splitting $A$, we make the length of the simple cycle in the down-dual graph including $A', A'', B$ even, and the length of the twisted cycles will become odd. Also if the incompatibility is happening in a branching $(N-1)$-simplex, we might need to split more than one of the cofaces of such simplex; this splitting of course will not change the number of branching and consequently will not affect the length of its corresponding cycle in the down-dual graph.      
   \end{proof}
 
   \begin{rem}
  Note that if splitting one/some of the $N$-simplex(es) divides its (their) free $(N-1)$-face(s) into two, then this will not have an effect on the other $N$-simplexes. But if we split the $N$-simplex in such a way that it divides its non-free face into two, due to the simplicial structure, this will also affect all other $N$-dimensional cofaces of the divided $(N-1)$ simplex, namely those which are lower adjacent to the original $N$-simplex (as shown in the next example). 
   \end{rem}

    \begin{ex}
In the  presented example in Figure 4, we show how to make simplicial complexes disorintable and the effect on the cycles of their corresponding down-dual graphs.

 We have two main rows of simplicial complexes (with blue) and the down-dual graph of each example is drawn at the bottom of each complex. The top row is the discrete  
Mobius strip and the bottom row corresponds to a tetrahedron.
  \begin{figure}[ht]
\centering
\includegraphics[width=0.65\textwidth]{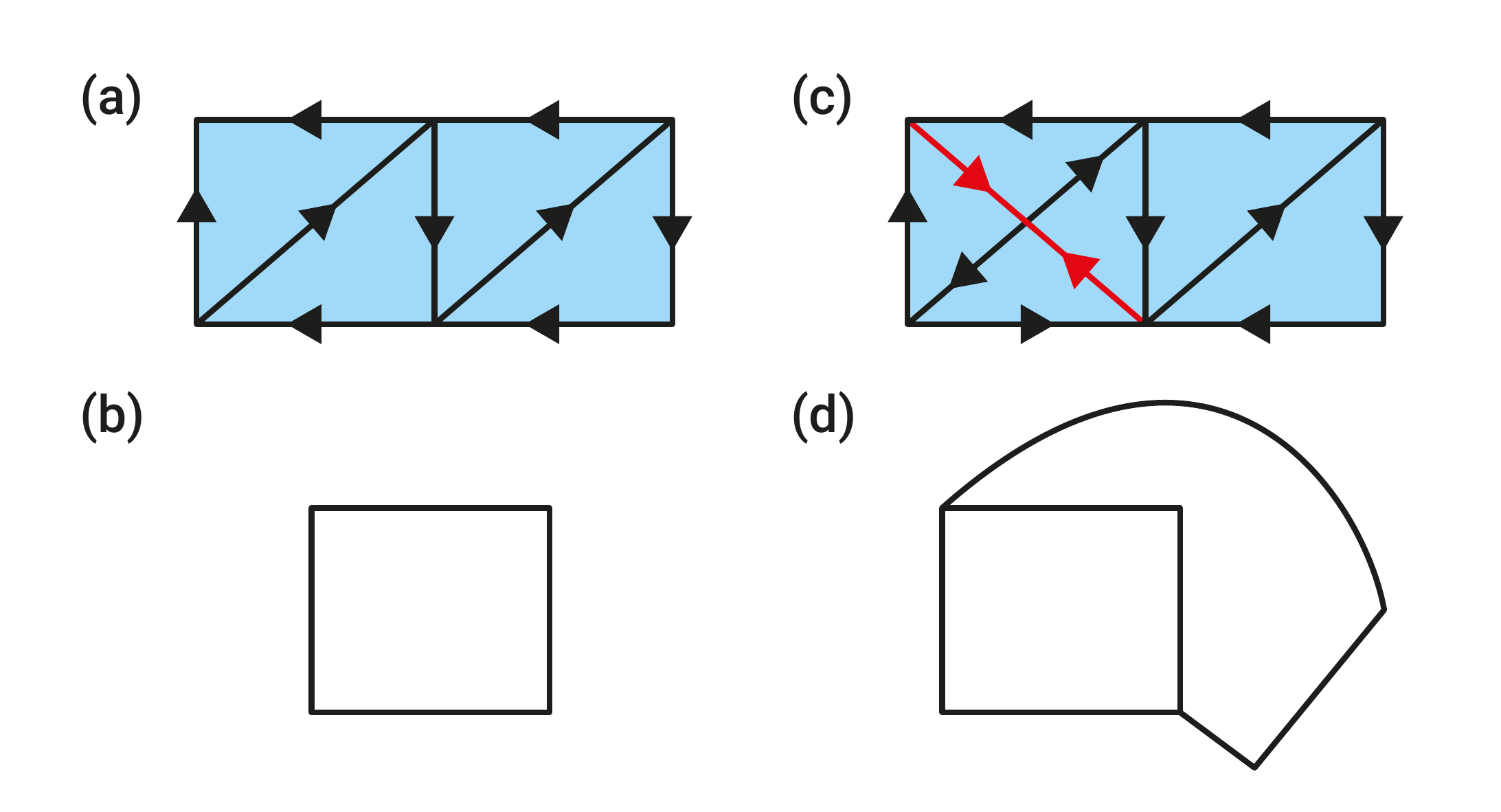}
\includegraphics[width=0.7\textwidth]{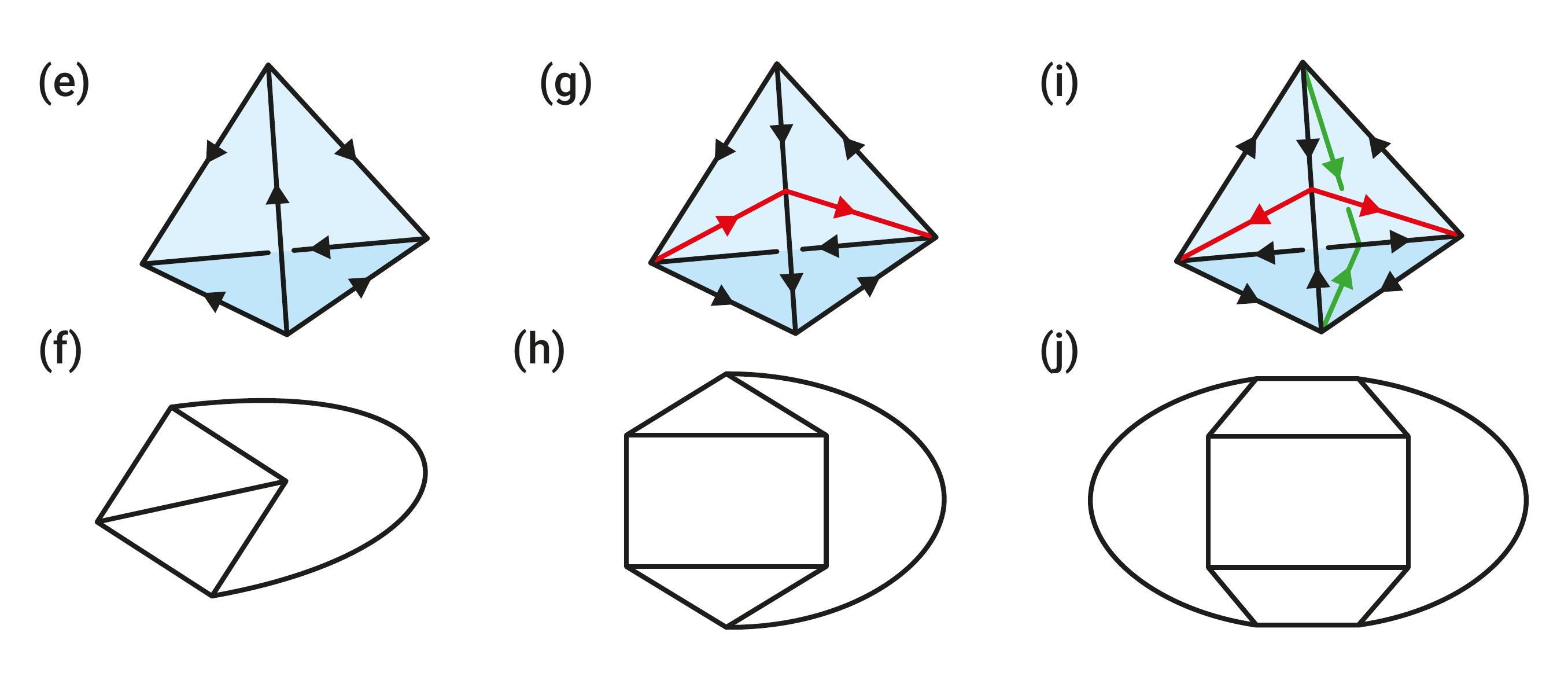}

\caption{Making simplicial complexes disorintable and the effect on the cycles of their corresponding down-dual graphs for Mobius strip and tetrahedron. 
The down-dual graph of each complex is drawn in its bottom. (a) is a non-disorientable 
Mobius strip as (b) has a twisting even cycle. (c) is obtained by splittings of two simplexes of (a), presented by red arrows, which is disorientable as in (d) twisting cycle(s) are even and simple cycle(s) are odd. Similarly, (e) and (g) are non-disorientable as based on respectively (f) and (h) they have odd simple cycles. (i) is obtained from (e) by splitting of all its four 2-simplexes and is disorientable as (j) has no odd simple simple.}
\end{figure}
   \end{ex}
\begin{lemma}

   The fewer odd cycles a graph has, the fewer splittings are required to make the graph bipartite; in other words, the closer the maximum eigenvalue is to 2.
\end{lemma}     
\begin{proof}
We should show that dividing an edge in an odd cycle of a graph into two by adding a vertex increases the maximum eigenvalue of the graph Laplacian. This can be demonstrated using the Rayleigh quotient, where an arbitrary real-valued function $f$ on the graph's vertex set is considered. If $O$ is an arbitrary odd cycle, there is always an edge $xy$ in the cycle that the function value at $x$ and $y$ has the same sign. By creating a new vertex $z$ (its degree is 2) between $x$ and $y$ with an opposite function sign, the Rayleigh quotient increases, as by adding $z$ the denominator  would increase by a factor of $2f^2(z)$  while the nominator increases by a factor of $(f(z)-f(x))^2+(f(z)-f(y))^2$ which is triviality bigger. Since this happens for any arbitrary real-valued function, the Rayleigh quotient as well as its supremum which is the maximum eigenvalue of the Laplacian increases and it becomes closer to the maximum possible value (i.e. two) and that is obtained by bipartite graphs with no odd cycles.    
\end{proof}
\begin{thm}
    In a $N$-simplicial complex without twisting cycles, 
   the fewer the number of (non-branching) simple odd cycles in its down dual graph, the fewer splittings of $N$- simplexes are required to make the complex disorientable; namely the closer is its maximum eigenvalue to the possible maximum eigenvalue (which  for the Normalized Laplacian is $N+1$).
\end{thm}
\begin{proof}
The proof is very similar to the one of the above Lemma. Here we use the down-dual graph (in dimension $N$) of the complex and we turn every simple (non-branching) odd cycle to even cycle one by one by adding one vertex as described before. Using Rayleigh quotient of the Laplacian in dimension $N-1$ and exact argument as above gives the result. As mentioned before, both of $L^{up}_{N-1}$ and $L^{down}_{N}$  have the same maximum obtained by the supremum of the Rayleigh quotient of functions defined on $N$-simplexes (i.e. vertices of the down dual graph).
\end{proof}
  
\textit{Conclusions} 

We have fully characterized disorientability of general simplicial complexes in terms of the parity of the length of cycles in their down dual graphs.
 As a direct result, every simplicial complex can become dissorientable by a finite number of splittings of some of the $N$-dimensional simplexes into two $N$-dimensional simplexes  in a similar manner that any graph can become bipartite by a finite number of splitting of some of its edges into two (i.e. making the odd edges even). Such splittings do not change the topology of the complex and the lower the number of splittings needed to make the complex disoriantable, the closer is the maximum eigenvalue of its Hodge Laplacian spectrum to the possible maximum eigenvalue.

We hope that this new perspective allows the extension of a range of theoretical and applied cycle-based methods from bipartite graphs to higher dimensional disorientable simplicial complexes. 

\bibliographystyle{unsrt}  
\bibliography{main}

\end{document}